\documentclass[11 pt]{amsart}
\usepackage[dvips]{graphicx}
\usepackage{amssymb}
\usepackage[english]{babel}
\usepackage{amsfonts,amsthm,amsmath}
\usepackage{epsfig}
\input{amssym.def}

\newtheorem{theorem}{Theorem}
\newtheorem{prop}[theorem]{Proposition}

\newtheorem{defi}[theorem]{Definition}
\newtheorem{remark}[theorem]{Remark}

\begin{document}

\title[An inverse spectral problem for surfaces.]{Determination of the Genus of Surfaces from the Spectrum of Schr\"odinger Operators attached to height functions}
\author{Brice Camus}
\address{Ludwig Maximilians Universit\"at M\"unchen,\\
Mathematisches Institut, Theresienstr. 39 D-80803
M\"unchen.\\Email: camus@math.lmu.de} \maketitle
\begin{abstract}
\noindent Using results on inverse spectral problems, in particular the so-called new wave invariants
attached to a classical equilibrium, we show that it is possible to determine the Morse index of height functions.
For compact Riemannian surfaces $M\subset \mathbb{R}^3$ this imply that we can retrieve the topology (via the genus).

Our results are independent from the choice of a metric on $M$ and
can be obtained from the choice of a 'generic' height-function.
For surfaces of genus zero, diffeomorphic to a 2-sphere, the
method allows to detect the convexity, or the local convexity of the surface.\medskip\\
keywords : Micro-local analysis; Schr\"odinger operators; Inverse
spectral problems.
\end{abstract}
\section{Introduction.}
\subsection{Basic definitions and setting.}
We are here interested in applications of the inverse spectral
problem for certain special differential operators in the
semi-classical regime. Let $X$ be compact (boundaryless)
Riemannian manifold equipped with a strictly positive density and
$\Delta_X$ the (positive) induced Laplace-Beltrami operator. In
particular, we consider the so called $h$-quantized Schr\"odinger
operator:
\begin{equation*}
P_{h}=h^2\Delta_X+V,\text{ on } L^2(X),
\end{equation*}
also called semi-classical Schr\"odinger operator, where the
potential $V$ is measurable and bounded from below on $X$.
\begin{remark}\rm{With some mild conditions on $V$, we could
also assume that $X$ is non-compact. But to simplify we stay in
the compact situation.}
\end{remark}
By a standard result, see \cite{Ber-Shu}, when $V$ is bounded from
below $P_{h}$ has a self-adjoint realization on a dense subset of
$L^2(X)$. To this quantum operator $P_{h}$ we can associate a
classical counterpart with the Hamiltonian, or total energy ,
function on the phase space:
\begin{equation*}
p(x,\xi)=||\xi||_x^2+V(x) \text{ on } T^* X.
\end{equation*}
Here the notation:
\begin{equation*}
||\xi||^2_x =\sum\limits_{i,j} g_{ij} \xi_i\xi_j,\, g_{ij}=G^{-1},
\end{equation*}
designs the norm (or scalar product at $x$) induced by the
Riemannian metric $G=g^{ij}(x)$ at $x$. We note $\Phi_{t}$ the
Hamiltonian flow of $H_{p}=\partial _{\xi}p.\partial_{x} -\partial
_{x}p.\partial_{\xi}$.

Here we are mainly interested in an asymptotic relation between
the semi-classical eigenvalues $\{\lambda_j(h): j\in\mathbb{N}\}$
of $P_h$:
\begin{equation*}
P_{h}\phi_j(x,h)=\lambda_j(h) \phi_j(x,h),\, \phi_j\in L^2(X),
\text{ as }h\rightarrow 0^+,
\end{equation*}
and the set of fixed point $\mathcal{P}$ (see below) for the map
$\Phi_{t}$ (viewed as a map on $\mathbb{R}_t\times T^*X$). We
refer to the introduction of \cite{Cam} for a general presentation
of this kind of relation between quantum and classical mechanics.

In the last section of this article we will work with compact
orientable surfaces $M$ of $\mathbb{R}^3$ and Schr\"odinger
operators $h^2\Delta_M+z$ attached
to height-functions $z$ on $M$.\medskip\\
\textbf{Spectral statistics.}\\
Consider an interval $I=[E_1,E_2]$ with $E_1<E_2$ and
$I(\varepsilon)=[E_1-\varepsilon,E_2+\varepsilon]$. For each
$\varepsilon>0$ the pullback $p^{-1}(I(\varepsilon))$ is a compact
subset of $T^*X$ and by a standard argument, see \cite{Ber-Shu},
it follows that the spectrum $\sigma (P_{h})\cap I(\varepsilon)$
is discrete and consists for each $h$ in a sequence:
\begin{equation*}
\lambda _{1}(h)\leq \lambda _{2}(h)\leq ...\leq \lambda _{j}(h),
\end{equation*}
of eigenvalues of finite multiplicities, if $\varepsilon$ and $h$
are positive. In general no formula is known to compute the
eigenvalues $\lambda_j(h)$ and to get pertinent information about
the spectrum (and the classical dynamics) it is interesting to
study the following spectral distributions:
\begin{equation}\label{Objet trace}
\Upsilon (E,h,\varphi)=\sum\limits_{\lambda _{j}(h)\in
I(\varepsilon)}\varphi (\frac{\lambda _{j}(h)-E}{h}),
\end{equation}
where $\varphi\in\mathcal{S}(\mathbb{R})$ is a test function, conveniently chosen, see below.

This kind of problem leads to a mathematically rigorous version of
the so-called Gutzwiller formula (see \cite{GUT}), a formula
intensively used in physics and quantum-chemistry. For example see
\cite{Haa} for various applications in physics and quantum chaos
with also many references. The so-called \textit{semi-classical
approximation} consists in studying the asymptotic behavior of
Eq.(\ref{Objet trace}) as the parameter $h$ tends to zero. For a
rigorous mathematical study of this problem a non-exhaustive list
of references is \cite{BU,Cdv,PU,Sj-Zw}.
\medskip\\
\textbf{Wave and new wave invariants.}\\
To study Eq.(\ref{Objet trace}), a classical approach (see \cite{BPU,BU} and section 3 of \cite{Cam}) is to study the asymptotic behavior, as $h\rightarrow 0^+$,
of the localized trace:
\begin{equation*}
\omega(E,h,t)=\mathrm{Tr}\, \left(\Theta(P_h) e^{-\frac{it}{h}
(P_h-E)}\right), \Theta\in C_0^\infty.
\end{equation*}
I will follow now the terminology used in \cite{Hez}. For $E$ regular it
is known since a while that $\Omega$ admits an asymptotic expansion of the
form:
\begin{equation*}
\omega(E,h,t)\sim \sum\limits_{j=-n}^{\infty} a_j(E,t) h^j, \text{
as } h\rightarrow 0^+.
\end{equation*}
\begin{defi}
The coefficients $a_j(E,t)$ are some distributions on the line
$\mathbb{R}_t$ and are called \textbf{wave invariants} of $P_h$.
\end{defi}
When $E=E_c$ is critical certain new coefficients appear in the
asymptotic expansion of $\omega(E,h,t)$ as $h$ tends to 0. In
particular, one can predict a general expansion in the form:
\begin{equation}\label{formal expand}
\omega(E_c,h,t)\sim
\sum\limits_{k=0}^{n-1}\sum\limits_{j=-n_0}^{\infty}
a_{j,k}(E_c,t) h^\frac{j}{p} \log(h)^k, \text{ as } h\rightarrow
0^+,
\end{equation}
for some $p\in \mathbb{N}^*$, see \cite{BPU} for the singularity
near $t=0$ and \cite{Cam} for many examples of new wave
invariants. In this work we only use the top-order coefficients of
Eq. (\ref{formal expand}) near a non-degenerate singularity (see
Theorem \ref{critical energy level} below).
\begin{defi}
These extra distributional coefficients appearing in Eq.(\ref{formal expand})
are called \textbf{new wave invariants}.
\end{defi}
In general, the top order coefficient w.r.t. $h$ of the
expansion involving the new wave invariants contains many information on
the shape of the symbol.
\section{Hypotheses and semi-classical results.}
We recall first the general result determining the wave invariants
at a regular energy level. Consider $X$ a closed smooth Riemannian
manifold, $n=\dim (X)$ and $V\in C^\infty(X)$ a positive
potential. Let $\{\lambda_j(h): j\in \mathbb{N}\}$ be the spectrum
of the Schr\"odinger operator $S_h=h^2 \Delta_X +V$ where
$\Delta_X$ is the Laplacian on $X$ (here given as a positive
operator). This spectrum of $S_h$ is always discrete, and with
finite multiplicities, when $X$ is compact or if $V$ is
'confining', i.e. $V(x)\rightarrow \infty$ when
$d(x,x_0)\rightarrow \infty$ for some $x_0\in X$. The
justification is that both conditions insure that the level-sets
$\Sigma_E$ defined below are compact and that the resolvent of
$S_h$ is a compact operator.

Let $\Phi_t$ be the Hamilton flow of $H(x,\xi)=||\xi||_x^2+V(x)$
on $T^* X$ and given $E>0$ define the energy surface:
\begin{equation*}
\Sigma_E =\{(x,\xi)\in T^*X: H(x,\xi)=E \}.
\end{equation*}
Recall that the flow can be viewed as a map
$\Phi_t:\Sigma_E\rightarrow \Sigma_E$ (conservation of the
energy). Also when the surfaces $\Sigma_E$ are compact the general
theory of differential equations insures that the flow is complete
(property of the maximal solutions of a Cauchy-problem). To
simplify notations we write $z=(x,\xi)\in T^*X$. We recall that
$E$ is regular (or non-critical) when $d H\neq 0$ everywhere on
$\Sigma_E$ and critical otherwise. Below, we use the subscript
$E_c$ to distinguish out critical values of $H$. The so-called
period manifold of $\Phi_t$ on $\Sigma_E$ is:
\begin{equation*}
\mathcal{P}= \{(T,z)\in\mathbb{R} \times \Sigma_E : \Phi_T(z)=z\}.
\end{equation*}
\textbf{At a non-critical energy level.}\\
When $E$ is non-critical, we have, see \cite{G-U} or \cite{BU},
the following general result concerning the wave invariants at the
energy $E$:
\begin{theorem}[\bf{Semi-classical trace formula at a regular level.}]\label{classic
expansion}$\,$\\
Assume that $E$ is regular and that the restriction of $\Phi_t$ to
$\Sigma_E$ is a clean flow (see section 2 of \cite{G-U}). Then
there exists a sequence of distributions on the real line,
$\{\gamma_k\}_k$, such that for every test function $\varphi$ with
Fourier-transform $\hat{\varphi}\in C_0^\infty$:
\begin{equation*}
\mathrm{Tr}\, \varphi (\frac{S_h-E}{h})= \sum\limits_{j=0}^\infty
\varphi( \frac{\lambda_j(h)-E}{h}) \sim \sum\limits_{j=1}^\infty
\gamma_j(\hat{\varphi}) h^{-n+j}c_j(h).
\end{equation*}
Moreover, the supports of the distributions $\gamma_j$ are contained in the \textbf{sets of periods} of the closed trajectories of $\Phi_t$ on $\Sigma_E$.
\end{theorem}
For a better description of the coefficients appearing in Theorem
\ref{classic expansion} we refer to \cite{G-U,BU} (we do not need
their explicit expressions here), see also \cite{D-G} for a
similar high-energy result concerning elliptic operators on a
compact manifold. Also, under our hypotheses the trace in Theorem
\ref{classic expansion} and the functional $\Upsilon(E,h,\varphi)$
are equal modulo a function of fast decay w.r.t. $h$. Such a
coefficient, of order $\mathcal{O}(h^\infty)$, is negligible in
semi-classical asymptotic expansions.

The idea we want to use here is that by a clever choice of $\mathrm{supp}(\hat{\varphi})$ we can eliminate the wave invariants appearing in Theorem \ref{classic expansion}:
\begin{itemize}
\item If $\hat{\varphi}$ is flat at the origin the set:
\begin{equation*}
\{ \{0\} \times \Sigma_E\} \subset \mathcal{P},
\end{equation*}
does not contribute.
\item If $\mathrm{supp}(\hat{\varphi})\subset [-T_0,T_0]$, for $T_0$ small enough then no periodic orbit:
\begin{equation*}
\{ (T,z): \Phi_T(z)=z\} \subset \mathcal{P},
\end{equation*}
will contribute to the asymptotic expansion.
\end{itemize}
\textbf{At a critical energy level.}
We allow now the presence of critical points for $H$ and we impose the type of singularity :\medskip\\
$(\mathcal{A}_{1})$\textit{ The potential } $V$ \textit{is a Morse function on} $X$.\medskip\\
A fortiori, in $I$ there is finitely many critical values
$E_c^1,...,E_c^l$ and in $p^{-1}(I)$ finitely many fixed points
$z_0^1,...,z_0^l$ of the energy function $p$.
\begin{remark}
\rm{The number of critical points $z_0^j=(x_0,0)$ is equal to the
number of critical energy levels. Otherwise $V$ would not be a
Morse function on $X$.}
\end{remark}
Next, we impose two conditions on our test function $\varphi$:\medskip\\
$(\mathcal{A}_2)$ \textit{$\hat{\varphi}$ is flat at 0, i.e. $\hat{\varphi}^{(j)}(0)=0$, $\forall j\in \mathbb{N}$}.\medskip\\
$(\mathcal{A}_3)$ \textit{For some sufficiently small $T$ we have $\mathrm{supp}(\hat{\varphi})\subset [-T,T]$.}\medskip\\
A fundamental property is that the singularity of
$\Upsilon(s,h,\varphi)$ as $s\rightarrow E_c$ describes partially
the singularity of $V$. In fact with conditions $(\mathcal{A}_2)$
and $(\mathcal{A}_3)$ we will only see the new wave invariants
attached to the critical point $z_0^j$ in $\Sigma_{E_j}$. We have:
\begin{theorem}[\bf{New wave invariants at a critical level}]\label{critical energy level}
$\,$\\
Under the conditions $(\mathcal{A}_{1})$, $(\mathcal{A}_{2})$ and $(\mathcal{A}_{3})$ we have:
\begin{equation*}
\mathrm{Tr}\, \varphi (\frac{S_h-E^j_c}{h})\sim
\sum\limits_{j=0}^\infty h^j c_j(\hat{\varphi}).
\end{equation*}
The leading coefficient is of the form:
\begin{equation*}
c_0(\varphi)=\frac{ e^{i\pi m_0/2 }}{(2\pi)} \int\limits_{\mathbb{R}} \frac{\hat{\varphi}(t)}{|\det (d\Phi_t(z_0^j) -\mathrm{Id})|^{\frac{1}{2}}},\, m_0\in \mathbb{Z}.
\end{equation*}
\end{theorem}
We refer to \cite{G-U,Cam,Cam0,KhD1} for a proof. Observe that $(\mathcal{A}_3)$ implicitly insures that $\det (d\Phi_t(z_0^j) -\mathrm{Id})\neq 0$.
Because of our implicit choice for $\varphi$, not all the new wave invariants are present in this formula.
The other new wave invariants are studied:
\begin{itemize}
\item In \cite{BPU}, near $t=0$.%
\item In \cite{Cam0,KhD1}, near a period of $d\Phi_t(z_0^j)$.
\end{itemize}
The explicit determination of all wave invariants, near a critical point of arbitrary signature, is a somehow complicated analytic problem involving
oscillatory integrals with degenerate phases. For an operator which is not a Schr\"odinger operator some new terms can generally appear at a period of $d\Phi_t(z_0)$ (see \cite{Cam0}).\medskip\\
\textbf{New wave invariants.} In our setting, the top-order
coefficient, given by the Duistermaat-Guillemin-Uribe density, is
indeed a smooth function as long as we stay away from any period
of the linearized flow at the point $z_0^j$. When $X=\mathbb{R}^n$
an explicit computation, done in \cite{KhD1} in a suitable system
of linear coordinates, shows that:
\begin{equation}\label{DGU}
d\nu_t(z_0)=\frac{1} {|\prod\limits_{j=1} ^r
\mathrm{sinh}(\alpha_j(z_0) t) \prod\limits_{j=r+1}^{n}
\sin(\alpha_j(z_0) t)|}.
\end{equation}
We must simply retain that the density $d\nu_t(z_0)$ determines:
\begin{itemize}
\item The signature $(n-r,r)$ of the Hessian of $V$ at $z_0$.
\item Eigenvalues $\alpha_j(z_0)$.
\end{itemize}
The last affirmation follows via Taylor-series and evaluation at several times.
\begin{remark}
\rm{In general, if the metric and the height function are unknown,
the spectral expectation determines only the numbers $\alpha_j$
and not the respective eigenvalues of $G(x_0)$ and $d^2V(x_0)$. A
similar indetermination is already valid for linear combinations
of harmonic oscillators on $\mathbb{R}^n$.}
\end{remark}
It follows that, when the potential is Morse-function on $X$, we
can retrieve the morse index of $X$ by several successive
applications of Theorem 6: we have only to cross finitely many
critical energy levels and to collect the index at each energy.
\section{Application to surfaces.}
There is a nice application to compact smooth surfaces $M\subset
\mathbb{R}^3$ equipped with a Riemannian metric (not necessarily
the metric of $\mathbb{R}^3$ restricted to $M$). We assume that
$M\subset \mathbb{R}^3$ is smooth, boundaryless, orientable and
that $M$ carries a smooth Riemannian metric $G$, fixed once for
all. We take $\Delta_M$ as the Laplace-Beltrami operator attached
to this metric (following the convention of geometers we may
assume that $\Delta_M$ is positive). Let us chose as potential $V$
a height function. We can assume $V$ to be positive, this can
always be achieved via a translation, $M$ being compact. If we
embed $M$ in $\mathbb{R}^3$, via some coordinates $(x,y,z)$, we
can chose $V$ as the projection on the $z$ axis. It is a standard
result of topology, see chapter 6 of \cite{Ban Hur}, that for
almost embedding $V$ will be a Morse function.

Then, for any choice of a smooth Riemannian metric on $M$, we have:
\begin{prop}\label{genus of surfaces}
Under the previous conditions on $V$, the semi-classical spectrum
of $P_h=h^2\Delta_M +V(x)$, defined as an unbounded operator on
$L^2(M)$, determines the topology of $M$.
\end{prop}
\begin{remark}\rm{Observe that the knowledge of the metric is not required. We only need a kinetic energy operator which is micro-locally elliptic and with a principal symbol
nowhere degenerated (see below). The knowledge of $V$ is also not required. We only need to recover the number of critical points of $V$ and their signature to conclude.}
\end{remark}
\noindent\textbf{Proof of Proposition \ref{genus of surfaces}.} We will use a variational argument w.r.t. the energy $E$.
Since $M$ is compact our potential has a maximum $E_{\mathrm{max}}$ and it will be sufficient to perform spectral estimates below $E_{\mathrm{max}}$.
Let $\lambda_j(h)$ be the spectrum of $P_h$, each eigenvalue being repeated according to it's multiplicity.

Since our potential is a Morse function, by Sard's theorem, we
obtain that the energy function $p(x,\xi)$ has only finitely many
critical values $E_c^j$, $j\in \{1,..., N\}$, attached to single
critical points. When $\mathrm{supp}(\hat{\varphi})$ is small
enough and does not contains the origin we have:
\begin{equation*}
\gamma(E,h,\varphi)\sim\left\{
\begin{matrix}
\mathcal{O}(h^\infty), \text{ for $E$ non-critical,}\\
c^j_0(\hat{\varphi}) +\mathcal{O}(h), \text{ for $E=E_c^j$ critical.}
\end{matrix}
\right.
\end{equation*}
Hence, the semi-classical spectrum determines each critical value $E_c^j$ of $V$.

Now for $j$ fixed we can use a simple micro-local argumentation.
The only critical point on $\Sigma_{E_c^j}$ is of the form
$z_0^j=(x_0^j,0)$ with $V(x_0^j)=E_c^j$. We pick a function in
$\psi\in C_0^\infty(T^* M)$ such that $0\leq \psi \leq 1$
everywhere and $\psi=1$ in a neighborhood of $z_0^j$. Always with
our conditions $(\mathcal{A}_2)$ and $(\mathcal{A}_3)$ on
$\mathrm{supp}(\hat{\varphi})$, we have:
\begin{equation*}
\Upsilon(E,h,\varphi)=\mathrm{Tr} \, \left( \psi^w(x,hD_x)
\varphi( \frac{P_h-E_c^j}{h})\right) +\mathcal{O}(h^\infty).
\end{equation*}
The important fact here is that on $\mathrm{supp}(1-\psi)$ there
is no critical point of $p$. Hence with condition
$(\mathcal{A}_2)$ and $(\mathcal{A}_3)$ we have:
\begin{equation*}
\mathrm{Tr} \, \left( (1-\psi^w(x,hD_x)) \varphi( \frac{P_h-E_c^j}{h})\right) =\mathcal{O}(h^\infty),
\end{equation*}
which easily follows from a non-stationary phase argument. On
$\mathrm{supp}(\psi)$, which can be chosen arbitrary small up to
an error of order $\mathcal{O}(h^\infty)$, we can use local
coordinates around $x_0^j$ and the Laplace operator has the form:
\begin{gather*}
-h^2 \sum\limits_{i,j} \sqrt{g} \frac{\partial}{\partial x_i} \frac{1}{\sqrt{g}}g_{ij} \frac{\partial}{\partial x_j}+V\\
= -h^2 \sum\limits_{i,j} g_{ij}(x) \frac{\partial}{\partial x_j}\frac{\partial}{\partial x_j}+V %
+h^2 \sum\limits_{i,j} \sqrt{g} \frac{\partial}{\partial x_i}(\frac{1}{\sqrt{g}}g_{ij})(x)\frac{\partial}{\partial x_j}.
\end{gather*}
Here we write the metric $G=g^{ij}$, $G^{-1}=g_{ij}$ and
$g=\det{G}$. Hence, in the sense of the $h$-calculus, we have
$p_h:=p_0+h p_1$ with a principal symbol:
\begin{gather*}
p_0(x,\xi)=\sum\limits_{i,j} g_{ij}(x)\xi_j\xi_i+V (x),
\end{gather*}
and a sub-principal symbol:
\begin{gather*}
p_1(x,\xi)= \sqrt{g}(x) \sum\limits_{i,j} \left(\frac{\partial}{\partial x_i} \frac{1}{\sqrt{g}}g_{ij}\right)_{|x} \, \xi_j.
\end{gather*}
Observe that $p_1=0$ at every point where $\xi=0$ and the
sub-principal symbol $p_1$ will play no essential r\^ole for the
estimates below. You can also use the convention that $h^2 \Delta$
is the quantization of $||\xi||^2$. This makes no difference for
the spectral estimates below.

Next, since $V$ is independent of the choice of the metric on $M$ we can freely assume
\footnote{Of course this choice is not isometric, in reality here should be the positive eigenvalues
of $G$ at $x_0^j$. But this choice is clearly sufficient to attain our objective. Observe that one could also locally trivialize the metric. This method is
a well known trick to derive Morse inequalities, e.g. for Witten-Laplacians.}  that near the origin:
\begin{equation*}
g_{ij}(x)=\mathrm{Id} +\mathcal{O}(||x||).
\end{equation*}
In this system of local coordinates we obtain that:
\begin{equation*}
\mathrm{Tr} \, \left( \psi^w(x,hD_x) \varphi( \frac{P_h-E_c^j}{h})\right)\sim C^j_0(\hat{\varphi}) +\mathcal{O}(h),
\end{equation*}
where:
\begin{equation*}
C^j_0(\hat{\varphi})= C \int\limits_{\mathbb{R}} \nu_t(z_0^j)
\hat{\varphi}(t)dt,\, C\in\mathbb{C^*}.
\end{equation*}
In our setting we can apply standard results on the linearized
flow, see \cite{A-M}, to compute the density $\nu_t(z_0^j)$. This
density, of the form given by Eq.(\ref{DGU}), determines the
number of positive and negative eigenvalues of $V$ at $x_0^j$.

In general position, define the Morse index of a critical point $x_0$ as the dimension of the negative eigenspace in $x_0$.
If $z$ is Morse function we denote by $N_j(z)$ the number of critical points of $z$ with index $j$. Since for a surface we have only $j=0,1,2$ we can retrieve the Euler characteristic:
\begin{equation*}
N_0(z)-N_1(z)+N_2(z)=\chi(M)=2(1-g(M)),
\end{equation*}
where $g(M)$ is the genus of $M$. $\hfill{\blacksquare}$
\begin{remark}\rm{Of course this approach is still valid in dimension $n>2$ but then the genus
is no more a sufficient topological invariant.
For a nice overview on Morse theory and indexes see \cite{Ban Hur} or \cite{Gra} for surfaces.
Observe that, for an unknown metric and an unknown height function, we can still retrieve
the Morse index of $V$ but not the Hessian of $V$. This is because in formula \ref{DGU} we can only retrieve the ratios $\alpha_j(z_0)$ of eigenvalues
of $G(x_0^j)$ and $d^2V(x_0^j)$ in a given system of coordinates. }
\end{remark}
The notion of Morse-Smale function (see, e.g., \cite{Ban Hur}
p.158) is here central. Morse-Smale functions are moreover dense
in every $C^r$-spaces ($r\geq1$) (Kupka-Smale-Theorem, p.159 and
remark 6.7 p.160 in \cite{Ban Hur}). To a Smale-Morse function is
attached the Morse-Smale-Floer complex and this complex is
isomorphic to the complex giving the singular homology (Theorem on
Morse homology, Theorem 7.4 of \cite{Ban Hur}). For a surface it
is known that the topology is given by the genus or the Euler
characteristic and this one is also given by the Euler-Poincar\'e
characteristic of the complex of homology.
That a certain Morse-function determines the topology of a surface is also contained in the book \cite{Gra} page 70.\medskip\\
\noindent\textbf{About the choice of a height function.}\\
For certain simple surfaces, e.g., convex surfaces diffeomorphic to a 2-sphere, there is no bad embedding since the height-function $z$ always shows up a strict minimum and a strict maximum. Such a
height function is also a \textit{perfect Morse function}, i.e. a Morse function with exactly 2 critical points attached respectively to a strict minimum and a strict maximum.
For a surface of genus 0 we can still have several critical points at the same critical value. This is not generic and unstable under a small perturbation of $z$.

For surfaces of higher genus some 'bad' embeddings are possible if
the height function is chosen transverse to a level set
(non-generic choice).
\begin{figure}[h!] \label{period functions}
\centering
   \begin{minipage}[l]{.36\linewidth}
       \epsfig{figure=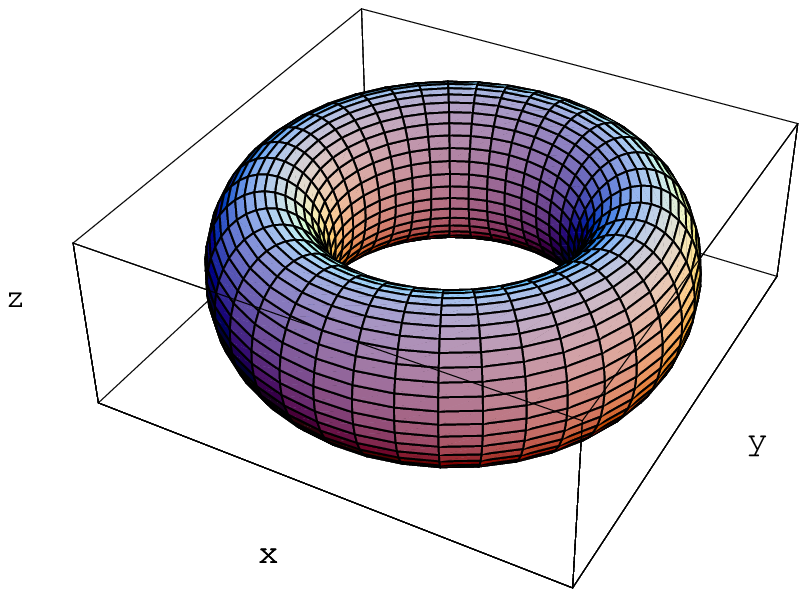,width=1. \textwidth}
   \end{minipage}
      \begin{minipage}[r]{.36\linewidth}
       \epsfig{figure=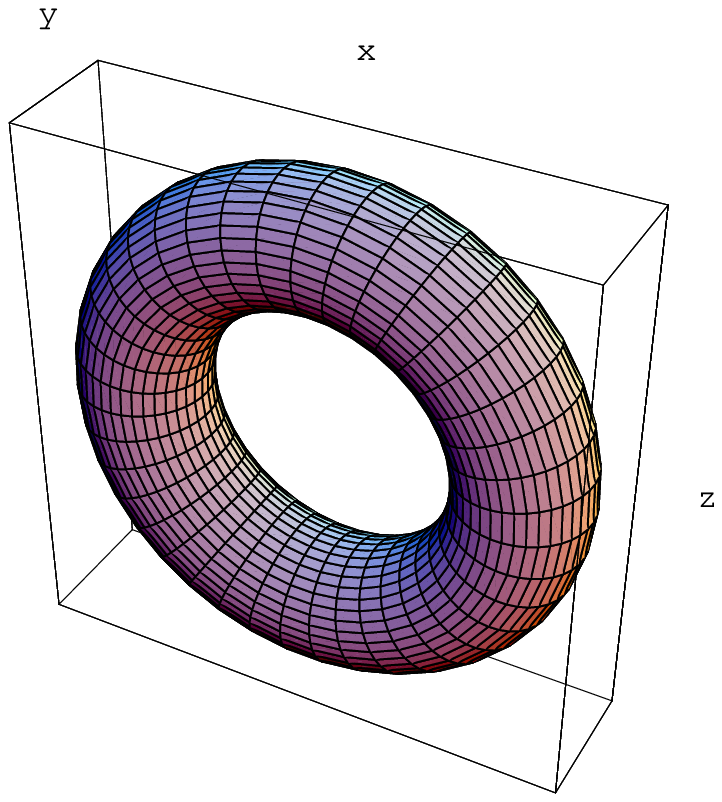,width=1. \textwidth}
   \end{minipage}
\caption{Embedding of a torus in $\mathbb{R}^3$.}
\end{figure}

The first embedding (here $z$ is the axe of symmetry) is not
favorable: the set of critical point consists of $2$ circles.
These circles are manifolds of critical points attached
respectively to a maxima and a minima of $V$ of energies
$E_\mathrm{min}<E_\mathrm{max}$. In that situation, a compact
manifold of critical points of dimension 1, we can here anyhow
apply the results of \cite{BPU} or\cite{KhD1}. If we still assume
that conditions $(\mathcal{A}_2)$ and $(\mathcal{A}_3)$ are
satisfied, each of these circles contributes as a 1-dimensional
submanifold in the 4-dimensional phase space:
\begin{equation*}
\Upsilon(E_\mathrm{min},\varphi,h) \sim C_1\int\limits_{\theta\in
\mathbb{S}^1}\int\limits_{t\in \mathrm{supp}(\hat{\varphi})}
\frac{h^{-\frac{1}{2}}\hat{\varphi}(t)}{|\sin(\alpha(\theta)t)|}
\frac{dt}{\sqrt{|t|}}d\theta+\mathcal{O}(h^{\frac{1}{2}}),\,
C_1\in\mathbb{C^*}
\end{equation*}
for the circle of minima and:
\begin{equation*}
\Upsilon(E_\mathrm{max},\varphi,h)\sim C_2\int\limits_{\theta\in
\mathbb{S}^1}\int\limits_{t\in \mathrm{supp}(\hat{\varphi})}
\frac{h^{-\frac{1}{2}}\hat{\varphi}(t)} {|\mathrm
{sinh}(\alpha(\theta)t)|}
\frac{dt}{\sqrt{|t|}}d\theta+\mathcal{O}(h^{\frac{1}{2}}),\,
C_2\in\mathbb{C^*}
\end{equation*}
for the circle of minima. Observe that the order w.r.t. $h$ is now $-1/2$.
\begin{remark}\rm{Both formulae for $\Upsilon(E_\mathrm{max},.)$
and $\Upsilon(E_\mathrm{min},.)$ easily follow from an application
of a stationary phase method with a compact manifold of critical
point and a non-degenerate transverse Hessian.}
\end{remark}
Here $\alpha(\theta)$ is the ratio of the eigenvalues of the
linearized operator in the transverse direction to
$\Sigma_{E_\mathrm{min}}\simeq \Sigma_{E_\mathrm{max}}\simeq
\mathbb{S}^1$ evaluated at the point $\theta\in\mathbb{S}^1$.
Observe that $\alpha(\theta)$ is negative for $E=E_\mathrm{max}$
(hyperbolic flow) and positive for $E=E_\mathrm{min}$ (periodic
flow). Observe that, only from the spectral estimates, we can
still see:
\begin{itemize}
\item Hyperbolic contributions: unstable equilibria at the maximal energy.
\item Trigonometric contributions: stable equilibria at the minimal energy.
\end{itemize}
The previous situation can be generalized for a smooth curve
$\gamma$, necessarily isomorphic to $\mathbb{S}^1$, of critical
points with a non-degenerate transverse Hessian at each point of
$\gamma$.
\begin{remark}
\rm{For a generic choice of the metric $G$ on $\mathbb{T}^2$ the
function $\theta \mapsto\alpha(\theta)$ is not constant along
$\mathbb{S}^1$. Unfortunately, the spectrum of the associated
Schr\"odinger operator $-h^2 \Delta_g+ z$ only determines the
average of the density along the circles. To get a better
description here requires to perform eigenfunction estimates. See,
e.g. \cite{BPU} for this point.}
\end{remark}

For the second embedding, where $x$ is the axe of symmetry, $z$ is a Morse function and we meet successively the critical points:
\begin{itemize}
\item a singularity of type $(0,2)$ : strict minimum,
\item a singularity of type $(1,1)$ : first saddle point,
\item a singularity of type $(1,1)$ : second saddle point,
\item a singularity of type $(2,0)$ : strict maximum.
\end{itemize}
This gives:
\begin{equation*}
\chi(\mathbb{T}^2)=1-2+1=0 \Rightarrow g(\mathbb{T}^2)=1.
\end{equation*}
The situation of the second embedding is generic and stable (e.g., w.r.t. a little deformation of the height function).
This allows to retrieve $\mathbb{T}^2$, up to a smooth deformation.
\begin{remark}\rm{Results concerning Morse-functions are not specific to height-functions and
Proposition \ref{genus of surfaces} can be generalized to
$h^2\Delta_M +V$ where $V$ is a Morse-function on $M$. The
interest here is the evident physical interpretation: the choice
of $V(x)=z(x)$, $x\in M$ is equivalent to put a particle, forced
to move freely on $M$ along geodesics, in a constant gravitation
field. Also a height-function gives a function independent of the
choice of the metric on $M$. For example, this is not the case of
a potential $V(x)$ depending (locally) on the geodesic distance
$d(x_0,x)$ on $M$. Such a potential depends on $G$ and can be
singular at conjugate points.}
\end{remark}
\noindent\textbf{Convexity and measure.}\\
Assume that $g(M)=0$ then if $M$ is convex every choice of a
height-function gives a perfect Morse function. But if $M$ is not
convex certain choice of the potential give locally a number of
critical points greater than 2, with the same Morse index.
\begin{figure}[h!]
\centering
   \begin{minipage}[l]{.50\linewidth}
       \epsfig{figure=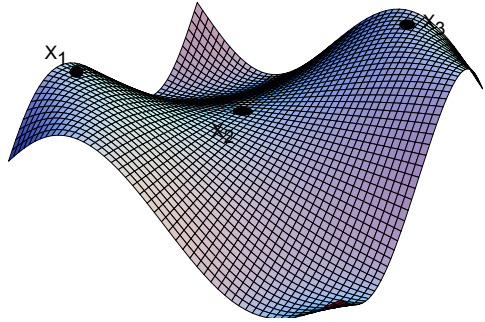,width=1. \textwidth}
   \end{minipage}
\caption{A non-convex surface: locally we see 3 critical points.}
\end{figure}

Observe that we still get 3 critical points by a small change for
$z$ and that we get an \textit{open condition}. In general to
detect the non-convexity of $M$ does not require global
spectral-estimates. Finally, for any gender $g(M)$ observe that
global spectral-estimates provide a lower bound for the euclidian
diameter of $M$ since we can retrieve the maximum and the minimum
of the restriction $z_{|M}$. This is just a lower bound as shows
figure 1.

From the point of view of statistical mechanics it could be
interesting to put a probability measure $\mu(z)$ on all choice
possible for $z$ and to average the spectral estimates with
respect to $\mu$. This is here simply equivalent to chose a
probability measure on $\mathbb{S}^2$ since the full problem is
invariant under translation. A similar construction is possible
for the choice of the metric $G$: if $G$ is in a \textit{bounded}
set of metrics $G_{\alpha}$, estimates given by conditions
$(\mathcal{A}_2)$ and $(\mathcal{A}_3)$ are still globally valid
and so are our conclusions. One could average
the results with respect to some probability measure $\mu(\alpha)$.\medskip\\
We could also obtain the contributions of a surface $M$ carrying a
flat section in the following sense:
\begin{quote} There exist an open subset $U\subset \mathbb{R}^3$ and a two-dimensional plane $P\subset \mathbb{R}^3$ such that
$U\cap M= U\cap P$.
\end{quote}
The flat section can be interpreted a 2-dimensional subset of
critical points when the height function $z$ is chosen
transversally to this section. The associated result is simply the
Lebesgue measure of the flat section. Observe that, as predicted
by the general theory of Morse-functions, this situation is not
generic and not stable under a small perturbation of $z$.\medskip\\
\textbf{Final Remark.} At a first look it might seem childish to
use semi-classical methods. But the 'high-energy' method (see e.g.
\cite{D-G}) is not working: when the energy $E$ is larger than the
maximum of the potential $E_{max}$ we have that the kinetic energy
is bounded from below by $||\xi||^2_x \geq E-E_{max}>0$. By
ellipticity of the Laplacian, it is not possible to produce any new wave invariant in this regime.\medskip\\
\textbf{Acknowledgments.} It is a great pleasure to thank George
Marinescu for useful discussions (and providing references)
concerning the Morse index. This work was partially supported by a
Deutsche Forschungsgemeinschaft Grant (D.F.G., the German research
foundation) \textit{Microlocal analysis applied to mathematical
physics and geometry}. The D.F.G. is greatly acknowledged for this
support.

\end{document}